\documentclass[11pt]{article}
\usepackage[hmargin=1in,vmargin=0.95in]{geometry}
\usepackage{amsmath,amssymb,amsthm,booktabs,array}
\usepackage{graphicx}
\usepackage{lmodern}
\usepackage{xcolor}
\definecolor{referenceblue}{RGB}{0,58,153}
\definecolor{abstractrule}{RGB}{172,181,192}
\usepackage[font=small,labelfont=bf,labelsep=period,justification=raggedright,singlelinecheck=false]{caption}
\usepackage[colorlinks=true,linkcolor=black,citecolor=black,urlcolor=referenceblue,bookmarksnumbered=true]{hyperref}
\hypersetup{pdfauthor={Feng Li and Pei Yu},pdftitle={Nine Limit Cycles in a Three-Dimensional Competitive Lotka--Volterra System}}

\newtheorem{theorem}{Theorem}
\newtheorem{proposition}[theorem]{Proposition}
\newtheorem{lemma}[theorem]{Lemma}

\newcommand{\tr}{\operatorname{tr}}
\newcommand{\Rea}{\operatorname{Re}}

\newcommand{\dd}{\mathrm d}
\newcommand{\HH}{\mathcal H}
\newcommand{\BB}{\mathcal B}
\title{Nine Limit Cycles in a Three-Dimensional Competitive Lotka--Volterra System}
\author{%
	\begin{tabular}{@{}c@{}}
	Feng Li\textsuperscript{1}\qquad  Yun Tian\textsuperscript{2} \qquad Pei Yu\textsuperscript{3}%
	\thanks{Pei Yu is the corresponding author (\nolinkurl{pyu@uwo.ca}).
		E-mail addresses:  \nolinkurl{lfeng@lyu.edu.cn} (F. Li), \nolinkurl{ytian22@shnu.edu.cn} (Y. Tian).}\\[0.45em]
	{\small\normalfont\textsuperscript{1} School of Mathematics and Statistics, Linyi University, Linyi, 276005, China}\\[0.12em]
	{\small\normalfont\textsuperscript{2} Department of Mathematics, Shanghai Normal University, Shanghai 200234, China}\\[0.12em]
	{\small\normalfont\textsuperscript{3} Department of Mathematics, Western University, London, N6A 5B7, Canada}
	\end{tabular}}
\makeatletter
\renewcommand{\@maketitle}{%
	\begin{center}
		{\LARGE\bfseries\@title\par}
		\vspace{0.75em}%
		{\large\@author\par}
	\end{center}%
}
\makeatother
\renewenvironment{abstract}{%
	\par\smallskip\noindent{\color{abstractrule}\rule{\linewidth}{0.4pt}}\par\smallskip
	\small\noindent\textbf{Abstract.}\enspace\ignorespaces
}{\par\medskip}
\date{}
\begin{document}
	\maketitle
	\begin{abstract}
		We construct a strictly competitive three-dimensional Lotka--Volterra system of Zeeman's class~27 with at least nine limit cycles. Eight bifurcate from the positive equilibrium. We compute the first eight focus values with the Tian--Yu recurrence, keeping the original interaction coefficients, and use a logarithmic identity to simplify each equation once the preceding ones have been solved. The first reduction gives a polynomial of degree seven. Exact rational and outward-rounded interval calculations verify an eighth-order fine focus and a nonsingular Jacobian for the eight unfolding coefficients. The local cyclicity is eight. For the same perturbation, the outermost small cycle and the boundary heteroclinic cycle both repel on the carrying simplex. Between them lies a ninth limit cycle.
	\end{abstract}
	
	{\small\noindent\textbf{Keywords.} Three-dimensional competitive Lotka--Volterra system, focus values, Hopf bifurcation, limit cycles, carrying simplex, heteroclinic cycle.\par}
	
	\section{Introduction}
	
	Can competing species coexist on a periodic orbit? The classical two species competitive Lotka-Volterra equations allow no positive periodic orbit, as follows from the Bendixson-Dulac criterion. With three species, competition can be cyclic. In the example of May and Leonard \cite{MayLeonard}, the identity of the dominant species changes repeatedly, and successive rounds take longer and longer. These oscillations arise as trajectories approach a cycle of connections between the single-species equilibria. The question of whether the three populations can follow an isolated periodic orbit takes us into the interior of the positive octant.
	
	Hirsch's carrying simplex gives a way to study both behaviors on the same surface \cite{Hirsch}. In three dimensions this invariant surface is homeomorphic to a closed triangle, and every positive periodic orbit lies on it \cite{Hirsch,XiaoLi}. Its boundary represents states in which at least one species is absent. Zeeman classified the robust nullcline configurations into 33 classes. Periodic orbits are possible only in classes~26-31 \cite{Zeeman,ZVD,JiangNiu}. On the carrying simplex of a class-27 system, the three single-species equilibria are saddles joined in a heteroclinic cycle. Thus there are two places to look for the dynamics that produce positive limit cycles. A Hopf bifurcation at the coexistence equilibrium can create cycles close to it. The flow near the boundary can help force a cycle farther out.
	
	Hofbauer and So, and later Xiao and Li, used this geometry to obtain one small cycle and a second cycle enclosing it \cite{HofbauerSo,XiaoLi}. Hofbauer and So conjectured that two was the upper bound. Lu and Luo reported three cycles in class~27 \cite{LuLuo27}, and Gyllenberg, Yan and Wang disproved the conjecture with three cycles in class~29, where there is no boundary heteroclinic cycle \cite{GYW}. In 2009, Gyllenberg and Yan returned to class~27 with a reported configuration of three small cycles and one enclosing cycle. They conjectured that a three-dimensional competitive Lotka--Volterra system could have at least five limit cycles \cite{GY}.
	
	The class-27 construction raises two questions about the same parameter values. As successive focus values are set to zero, how many small cycles can bifurcate from the positive equilibrium? What happens to the boundary stability at those parameters? In a separate study, Gyllenberg and Yan found a class-27 Hopf point at which the boundary characteristic number and the first two focus values vanish, while the third focus value is nonzero \cite[Theorem~2.1]{GYDEA}. At that point the sign criterion for boundary stability is inconclusive.
	
	Computing enough focus values became a central part of the problem. Tian and Yu developed a recurrence that computes the center manifold and the Hopf normal form together. In their Lotka--Volterra application they discussed a possible configuration with seven small cycles and one enclosing cycle, conditional on finding a seventh-order fine focus and obtaining the required boundary sign \cite[Sect.~3.2]{TianYu}. Yu, Han and Xiao subsequently re-examined the earlier class-27 examples. They found a Lyapunov sign error in Lu and Luo's construction and a problem with the enclosing-cycle argument for the first family of Gyllenberg and Yan. At a different critical point in the second family of Gyllenberg and Yan, they established three small cycles and one enclosing cycle. They also obtained four small cycles in separate class-27 examples \cite[Sects.~3.5, 3.7 and~4.1]{YHX}.
	
	Murakami later gave an example with two small cycles and one enclosing cycle \cite{Murakami}. Alongside the work on Lyapunov coefficients and center conditions \cite{WangHuang,WLHS,JLWH}, Li and Jiang established three-cycle lower bounds in classes~27--30 \cite{LJ}. The four-cycle constructions now include class~29 and, in a preprint of Hu, Lu and Luo, class~28 \cite{WHW,HLL}. More recently, in a preprint, Wang, Shao and Dai gave rational competitive systems with five, six, seven and eight small-amplitude hyperbolic limit cycles \cite[Theorem~1]{WSD}. We ask whether a class-27 system can have eight small cycles and a cycle that surrounds them all.
	
	We seek a perturbation for which the outermost small cycle and the boundary heteroclinic cycle both repel on the carrying simplex. When the two signs hold for the same perturbed system, we can trap a trajectory in the region between them and look for an additional periodic orbit there.
	
	We place the positive equilibrium at \(E=(1,1,1)\) and write
	\begin{equation}\label{eq:model}
		\dot x_i=x_i\sum_{j=1}^3 A_{ij}(x_j-1),\qquad i=1,2,3,\qquad A_{ij}<0 .
	\end{equation}
	We use Tian and Yu's recurrence to compute the focus values and keep the entries of \(A\) as parameters. After restricting to the zero set of the preceding focus values, a logarithmic identity simplifies the next focus equation. Restricted derivatives then allow us to solve for one matrix entry at a time. We verify a Hopf point with \(v_1=\cdots=v_7=0\) and \(v_8>0\), together with a nonsingular Jacobian for the unfolding. Finally, we show that the repelling boundary persists under the perturbation that gives the eight small cycles and makes the outermost one repelling. This yields the following result.
	
	\begin{theorem}\label{thm:main}
		There exists a strictly competitive three-dimensional Lotka--Volterra system of Zeeman class~27 with at least nine distinct limit cycles in the positive octant. Eight small-amplitude limit cycles bifurcate from the positive equilibrium, which is an eighth-order fine focus at the critical parameter values. An additional limit cycle encloses all eight on the carrying simplex.
	\end{theorem}
	
	\section{The Hopf condition and computation of the focus values}\label{sec:recurrence}
	
	For system \eqref{eq:model}, write
\[
A=\begin{pmatrix}a&b&c\\ d&e&f\\ g&h&k\end{pmatrix},
\qquad \tau=\tr A,\qquad s=\sigma_2(A),\qquad
D=\det A-\tau s ,
\]
where \(\sigma_2\) is the sum of the principal minors of order two. 
Then the characteristic equation of $A$ can be rewritten as $(\lambda-\tau)(\lambda^2+s)-D=0$.
Solving \(D=0\) yields
\begin{equation}\label{eq:hopf}
	f=\frac{(a+e)(a+k)(e+k)-bd(a+e)-cg(a+k)-cdh}
	{bg+(e+k)h}.
\end{equation}
Indeed, \(D\) is affine in \(f\), with coefficient \(\delta=bg+(e+k)h>0\) for a strictly negative matrix. The condition \(s>0\) gives eigenvalues \(\tau,\eta,-\eta\), with \(\eta=i\sqrt{s}\) and \(\tau<0\).

\subsection{The first focus value}
 We fix \(k=-1\) by a positive rescaling of time. 
Put \(X=x-\mathbf 1\), and write \(\odot\) for the componentwise product. Then \(\dot X=AX+X\odot AX\). Define
\[
r(\eta)=
\begin{pmatrix}
	bf-c(e-\eta)\\ cd-f(a-\eta)\\ (a-\eta)(e-\eta)-bd
\end{pmatrix},
\quad q=\frac{r(\eta)}{Q(\eta)}.
\]
 where \(Q(\eta)=r_1(\eta)\) is the first entry of \(r(\eta)\).
The imaginary part of \(Q(\eta)\) is \(c\sqrt{s}\ne0\). Let \(e_j\) denote the \(j\)th coordinate vector, and put
\[
K_\tau=\frac{A^2+sI}{\tau^2+s},\qquad
q_\tau=\frac{K_\tau e_1}{e_1^{\mathsf T}K_\tau e_1},
\qquad T=(q,\bar q,q_\tau).
\]
The denominator in \(q_\tau\) is positive since \(A\) is strictly negative and \(s>0\). The columns of \(T\) correspond to the distinct eigenvalues \(\lambda_1=\eta,\lambda_2=-\eta,\lambda_3=\tau\). Let \(\ell\) be the first row of \(T^{-1}\).

With the normalization given below, the first focus value is
\begin{equation}\label{eq:v1}
	v_1=\Rea\left\{\frac{\eta(\tau-\eta)}{2\eta-\tau}
	\ell\big[\bar q\odot K_\tau(q\odot q)\big]\right\}.
\end{equation}
To obtain its numerator in the original parameters, reduce the following expressions modulo \(\eta^2+s\):
\begin{align*}
	U(\eta)&=e_1^{\mathsf T}(A+\eta I)(A-\tau I)
	\big[r(-\eta)\odot (A^2+sI)(r(\eta)\odot r(\eta))\big]
	= U_0+U_1\eta,\\
	W(\eta)&=(2\eta-\tau)Q(\eta)=W_0+W_1\eta .
\end{align*}
Using \(\ell=e_1^{\mathsf T}(A+\eta I)(A-\tau I)/[2\eta(\eta-\tau)]\) and rationalizing \eqref{eq:v1} gives
\begin{equation}\label{eq:num}
	v_1=-\frac{U_0W_0+sU_1W_1}
	{2(\tau^2+s)(\tau^2+4s)\,[Q(\eta)Q(-\eta)]^2}.
\end{equation}
The denominator is strictly positive, and the expression remains rational after substituting \eqref{eq:hopf}.
	
\subsection{The Tian--Yu Hopf normal-form recurrence}

\label{sec:tianyu}
We use the Tian--Yu recurrence to compute the center manifold and normal form together \cite[Eq.~(7) and Theorem~1]{TianYu}. In the coordinates \(y=T^{-1}X\), introduce the formal parameterization
\[
y=Y(z,w)=\sum_{m,n\ge0}Y_{mn}z^mw^n,\qquad
w=\bar z,\quad Y_{10}=e_1,\quad Y_{01}=e_2,\quad Y_{00}=0.
\]
We treat \(z,w\) as independent in the coefficient computations and impose \(w=\bar z\) for the real system. The first two components give the center coordinates as functions of \(z,w\), and the third gives the stable coordinate on the manifold. Write \(Y_{mn,i}\) for component \(i\) of \(Y_{mn}\). The reality condition is \(Y_{nm}=(\overline{Y_{mn,2}},\overline{Y_{mn,1}},\overline{Y_{mn,3}})^{\mathsf T}\). In the translated variables \(X=H(z,w)=TY(z,w)\), the coefficients are \(H_{mn}=TY_{mn}\). We set coefficients with negative indices to zero and denote the reduced field \((\dot z,\dot w)^{\mathsf T}\) by \(\mathcal G=(\mathcal G_1,\mathcal G_2)^{\mathsf T}\), where
\[
\mathcal{G}_1=\eta z+\sum_{j\ge1}c_jz^{j+1}w^j,\qquad
\mathcal{G}_2=-\eta w+\sum_{j\ge1}\bar c_jz^jw^{j+1},
\qquad v_j=\Rea c_j .
\]
Set \(\Lambda=\operatorname{diag}(\eta,-\eta,\tau)\) and \(\mathcal Q(U,V)=T^{-1}\big[(TU)\odot A(TV)\big]\). Comparison of coefficients in the invariance equation \((D_{(z,w)}Y)\mathcal G=\Lambda Y+\mathcal Q(Y,Y)\) gives, for \(m+n\ge2\), the forcing term
\begin{align}
	R_{mn}={}&
	\sum_{\substack{0\le u\le m,\ 0\le v\le n\\0<u+v<m+n}}
	\mathcal Q(Y_{uv},Y_{m-u,n-v})\notag\\
	&-\sum_{1\le j\le\lfloor(m+n-2)/2\rfloor}
	\big((m-j)c_j+(n-j)\bar c_j\big)Y_{m-j,n-j}.\label{eq:R}
\end{align}
Let \(R_{mn,i}\) be component \(i\) of \(R_{mn}\). The homological equation for component \(i\) has the factor \((m-n)\eta-\lambda_i\). At resonance we set the transformation coefficient to zero and retain the forcing as a normal-form coefficient. This gives
\begin{equation}\label{eq:H}
	Y_{mn,i}=
	\begin{cases}
		R_{mn,i}/\big((m-n)\eta-\lambda_i\big), & (m-n)\eta\ne\lambda_i,\\
		0, & (m-n)\eta=\lambda_i,
	\end{cases}
	\qquad c_j=R_{j+1,j,1}.
\end{equation}
In polar coordinates \(z=re^{i\theta}\), the formal radial equation on the Hopf chart is \(\dot r=r\sum_{j\ge1}v_jr^{2j}\). We therefore carry the calculation through degree~17 to obtain \(v_1,\ldots,v_8\).

To retain the original parameters in the successive eliminations, we perform the coefficient arithmetic over a rational function field. Substitute \(k=-1\) and \(f=f_{\rm H}(a,b,c,d,e,g,h)\) from \eqref{eq:hopf}, and represent each scalar \(x=x_0+\eta x_1\) by \((x_0,x_1)\). Multiplication and scalar inversion then take the forms
\[
(x_0,x_1)(y_0,y_1)=(x_0y_0-sx_1y_1,\;x_0y_1+x_1y_0),\qquad
(x_0,x_1)^{-1}=\frac{(x_0,-x_1)}{x_0^2+sx_1^2}.
\]
Complex conjugation changes the sign of the second component. These rules, applied to \eqref{eq:R}--\eqref{eq:H}, give \(v_j\) as the first component of \(c_j\). All operations remain in the rational function field in \(a,b,c,d,e,g,h\), and every inverse required by the recurrence exists on the admissible Hopf chart.

We can now recover \eqref{eq:v1} directly from the recurrence. Put \(u=T^{-1}(q\odot q)\). At quadratic order, \(Y_{11}=0\) and
\[
Y_{20}=\left(u_1,\frac{u_2}{3},\frac{\eta u_3}{2\eta-\tau}\right)^{\mathsf T}.
\]
The cubic resonance yields \(c_1=\eta\ell[\bar q\odot(H_{20}-q\odot q)]\). The first component of \(Y_{20}-u\) vanishes. The second contributes the purely imaginary term \(-2\eta|\bar\ell(q\odot q)|^2/3\) to \(c_1\). The third gives \eqref{eq:v1}, since \(u_3q_\tau=K_\tau(q\odot q)\).

\section{Successive elimination and center conditions}\label{sec:elimination}

\subsection{Successive reduction of the eight focus values}

The logarithmic structure of the LV system simplifies the higher-order calculations after each focus equation is imposed. Let \(\HH_j=\{v_1=\cdots=v_j=0\}\) on the Hopf chart, and let \(\HH_0\) be the entire chart. Put \(L(z,w)=\log(\mathbf 1+H(z,w))\) componentwise, with \(L_{mn}\) denoting its vector coefficient of \(z^mw^n\). Then \eqref{eq:model} gives \((D_{(z,w)}L)\mathcal G=AH\). The coefficient of \(z^nw^n\), for \(n\ge1\), gives
\begin{equation}\label{eq:log}
	AH_{nn}=2\sum_{m=1}^{n-1}(n-m)v_mL_{n-m,n-m}.
\end{equation}
In particular, we obtain
\[
\begin{aligned}
	H_{11}&=0,\qquad L_{11}=-q\odot\bar q,\qquad H_{22}=2v_1A^{-1}L_{11},\\
	H_{33}&=2A^{-1}(2v_1L_{22}+v_2L_{11}),\\
	H_{44}&=2A^{-1}(3v_1L_{33}+2v_2L_{22}+v_3L_{11}).
\end{aligned}
\]
The eigenvalues of \(A\) are all nonzero on the Hopf chart, so \(H_{nn}=0\) on \(\HH_{n-1}\). Hence, at order \(j\), the preceding focus equations remove \(H_{11},\ldots,H_{jj}\). The preceding \(c_m\) are then purely imaginary. Their imaginary parts must still be retained in the recurrence.

Write \(\BB(U,V)=U\odot AV+V\odot AU\). At the resonance \((j+1,j)\), the first component of the second sum in \eqref{eq:R} vanishes, since \(Y_{j+1-m,j-m,1}=0\) by \eqref{eq:H}. In the first sum, we combine each product with its reverse order. As the first row of \(T^{-1}\) is \(\ell\), this yields
\begin{equation}\label{eq:reduced}
	\begin{aligned}
		T_{jm}&=\sum_{\substack{0\le u\le m\\u\ne\lceil m/2\rceil}}
		\BB(H_{u,m-u},H_{j+1-u,j-m+u}),\qquad 1\le m\le j,\\
		S_j&=\sum_{m=1}^j T_{jm},\qquad
		D_j=\sum_{n=1}^j\BB(H_{nn},H_{j+1-n,j-n}),\\
		v_j&=\Rea\bigl(\ell(S_j+D_j)\bigr),\qquad
		\widehat v_j=\Rea(\ell S_j).
	\end{aligned}
\end{equation}
Here \(m\) is the smaller total degree of the two factors. The pair indexed by \(u=\lceil m/2\rceil\) is the only one in each group that contains a diagonal coefficient, and it is included in \(D_j\). Table~\ref{tab:focusforms} gives \(S_j\) at each order, and Appendix~\ref{app:terms} lists the individual products. The coefficients \(H_{mn}=TY_{mn}\) used here are computed from the full recurrence \eqref{eq:R}--\eqref{eq:H}. On \(\HH_{j-1}\), the terms in \(D_j\) vanish, so \(v_j=\widehat v_j\) there.

\begin{table}[htbp]
	\caption{Factored focus-value expressions in the original matrix entries. At order \(j\), \(D_j\) consists of \(j\) diagonal products given by \eqref{eq:reduced}. Their coefficients are determined by \eqref{eq:log}.}
	\label{tab:focusforms}
	\centering
	\small
	\renewcommand{\arraystretch}{1.08}
	\begin{tabular}{@{}cclc@{}}
		\toprule
		\(j\)&Degree&\(S_j\) in \(v_j=\Rea(\ell(S_j+D_j))\)&Terms in \(S_j\)\\
		\midrule
		1&3 &\(T_{11}\)&1\\
		2&5 &\(T_{21}+T_{22}\)&3\\
		3&7 &\(T_{31}+T_{32}+T_{33}\)&6\\
		4&9 &\(T_{41}+T_{42}+T_{43}+T_{44}\)&10\\
		5&11&\(T_{51}+T_{52}+T_{53}+T_{54}+T_{55}\)&15\\
		6&13&\(T_{61}+T_{62}+T_{63}+T_{64}+T_{65}+T_{66}\)&21\\
		7&15&\(T_{71}+T_{72}+T_{73}+T_{74}+T_{75}+T_{76}+T_{77}\)&28\\
		8&17&\(T_{81}+T_{82}+T_{83}+T_{84}+T_{85}+T_{86}+T_{87}+T_{88}\)&36\\
		\bottomrule
	\end{tabular}
\end{table}

\paragraph{Orders one and two.}
At the first order, \(v_1=\Rea\ell\BB(\bar q,H_{20})\), which gives the rational expression \eqref{eq:num}. At the second order, the full value and its restriction are
\[
\begin{aligned}
	S_2={}&\BB(\bar q,H_{31})+\BB(H_{02},H_{30})+\BB(H_{20},H_{12}),\\
	v_2={}&\Rea\ell\bigl[S_2+\BB(q,H_{22})\bigr],\\
	\widehat v_2={}&\Rea\ell S_2\qquad\text{on }v_1=0.
\end{aligned}
\]
To evaluate these terms explicitly, let \(\mathcal L_l\) be the diagonal matrix of the three multipliers in \eqref{eq:H} when the difference of the Taylor indices is \(l\). We obtain
\[
\begin{aligned}
	Y_{20}&=\eta\mathcal L_2T^{-1}(q\odot q), & H_{02}&=\bar H_{20},\\
	c_1&=\ell\BB(\bar q,H_{20}), &
	Y_{21}&=\mathcal L_1T^{-1}\BB(\bar q,H_{20}),\\
	Y_{30}&=\mathcal L_3T^{-1}\BB(q,H_{20}), &H_{12}&=\bar H_{21},\\
	Y_{31}&=\mathcal L_2\bigl[T^{-1}\bigl(\BB(q,H_{21})+\BB(\bar q,H_{30})\bigr)-2c_1Y_{20}\bigr].
\end{aligned}
\]
When \(v_1=0\), we have \(H_{22}=0\), and the term \(-2c_1Y_{20}\) remains in \(Y_{31}\).

\paragraph{Order three.}
On the full Hopf chart,
\[
\begin{aligned}
	S_3={}&\BB(H_{01},H_{42})+\BB(H_{02},H_{41})+\BB(H_{20},H_{23})\\
	&+\BB(H_{03},H_{40})+\BB(H_{12},H_{31})+\BB(H_{30},H_{13}),\\
	D_3={}&\BB(H_{11},H_{32})+\BB(H_{22},H_{21})+\BB(H_{33},H_{10}).
\end{aligned}
\]
Imposing \(v_1=v_2=0\) gives \(H_{22}=H_{33}=0\) and hence \(v_3=\widehat v_3\). At the fourth order, the same identity removes \(H_{44}\) after \(v_3=0\). Continuing in this way gives the reduction through \(H_{88}\).

Formula \eqref{eq:reduced} defines \(\widehat v_j\) on the full chart using the same \(H_{mn}\) as for \(v_j\). By \eqref{eq:log}, each \(H_{nn}\) in \(D_j\) is a linear combination of \(v_1,\ldots,v_{n-1}\) with analytic coefficients. Consequently \(\widehat v_j-v_j=-\Rea(\ell D_j)\) belongs to the ideal \((v_1,\ldots,v_{j-1})\) of analytic functions. At a common zero, differentiating this identity gives a unit lower triangular row transformation between the Jacobians of the first seven full and reduced values.

\subsection{Successive solution of the focus equations}

Recall \(p=(p_1,\ldots,p_7)=(a,b,c,d,e,g,h)\). Throughout this subsection,
derivatives with respect to \(p\) are taken on this chart, and
\(f_{\rm H}\) is recomputed after each substitution. Section~\ref{sec:validation} uses
interval contraction to certify the exact common zero \(p_*\) of
\(F=(v_1,\ldots,v_7)\).

\paragraph{Hopf equation solves \(f=A_{23}\).}
With the other entries fixed, we have \(D=\delta(f-f_{\rm H}(p))\), where
\(\delta=bg+(e-1)h>0\) does not depend on \(f\). The Hopf equation
\(D=0\) is therefore affine in \(f\), with the unique solution
\eqref{eq:hopf}.

\paragraph{Evaluation after substitution.}
Write \(t_j=(p_j,\ldots,p_7)\) for the entries before step \(j\), with
\(t_8\) the empty tuple. Start with \(E_0(p)=p\). Once the \(j\)th
equation has been solved as \(p_j=\phi_j(t_{j+1})\), define
\[
E_j(t_{j+1})=E_{j-1}(\phi_j(t_{j+1}),t_{j+1}).
\]
For \(1\le j\le7\), the map \(E_j\) records all seven entries after the first \(j\)
solutions have been substituted. At step \(j\), we solve \(G_j(t_j)=0\), where
\begin{equation}\label{eq:gj}
	G_j(t_j)=v_j(E_{j-1}(t_j))
	=\widehat v_j(E_{j-1}(t_j)).
\end{equation}
The preceding values \(v_1,\ldots,v_{j-1}\) vanish at the substituted
point, so the second equality follows from \eqref{eq:reduced}.
The matrix \(A\), the coefficients \(H_{uv}\), and the row vector
\(\ell\) are evaluated there as well.

\paragraph{The restricted derivative and local solvability.}
Lemma~\ref{lem:cert} first establishes an exact common zero of the full
focus-value functions in the interior of the box. We then construct the
graphs \(E_j\) successively through this point. On each preceding zero
locus, omitting the diagonal terms leaves the focus value unchanged.
The derivatives of the full functions on the Hopf chart determine the
tangent to that locus.

For an entry \(u\) of \(p\), write \(v_{i,u}=\partial_u v_i\) and let
\(J_j\) be the leading \(j\)-by-\(j\) block of the Jacobian \(DF\), with
\(\det J_0=1\). The notation \((E_{j-1})_{1:j-1}\) denotes the first
\(j-1\) components of \(E_{j-1}\).
Differentiate the first \(j-1\) equations with respect to \(p_j\),
holding the trailing arguments fixed. The derivatives of the entries
already solved for satisfy
\[
J_{j-1}\frac{\partial(E_{j-1})_{1:j-1}}{\partial p_j}
=-\frac{\partial(v_1,\ldots,v_{j-1})}{\partial p_j}.
\]
The chain rule then gives
\begin{equation}\label{eq:schur}
	\partial_{p_j}G_j
	= v_{j,p_j}
	-\frac{\partial v_j}{\partial(p_1,\ldots,p_{j-1})}\,
	J_{j-1}^{-1}
	\frac{\partial(v_1,\ldots,v_{j-1})}{\partial p_j}
	=\frac{\det J_j}{\det J_{j-1}}.
\end{equation}
All derivatives and Jacobian blocks on the right are evaluated at \(E_{j-1}(t_j)\).
For \(j=1\), the correction term is absent. The ratios are precisely the
unpermuted Gaussian pivots of \(DF\), which the integer dyadic calculation
encloses over the whole rational box. At each step, the product of the
first \(j-1\) pivots certifies \(\det J_{j-1}\ne0\). The eigenbasis,
scalar homological, and Hopf-chart denominators do not vanish near the
box, so the full focus values are analytic there. The trailing entries
of \(p_*\) satisfy \(G_j=0\), and the nonzero restricted derivative allows
us to apply the implicit-function theorem. This gives the next analytic
solution graph on a small enough neighborhood that it stays inside the
certified box. At the last step it gives an isolated simple root. Notice that
\(v_{j,p_j}\) is only the first term of the derivative along the preceding
solution graphs in \eqref{eq:schur}.

\paragraph{Step 1. \(v_1=0\) solves \(a=A_{11}\).}
We evaluate the numerator \(N_1=-(U_0W_0+sU_1W_1)\) of \eqref{eq:num} at
\(A(p)\). Because \(f-f_{\rm H}\) is monic and linear in \(f\), eliminating
\(f\) by a resultant amounts to the same substitution. Cancellation gives
the septimic equation \(S_7^{\rm alg}(a,b,c,d,e,g,h)=0\), derived in
Section~\ref{sec:septimic}. Its leading coefficient in \(a\) is \(-gh(e-1)^3>0\).
The denominator in \eqref{eq:septimic} is nonzero, and the enclosure
\(\partial_aG_1\in[1.0162706,1.0162707]\times10^{-2}\) shows that the root
at \(p_*\) is simple. We may therefore write \(a=\phi_1(b,c,d,e,g,h)\) and set
\[
E_1(b,c,d,e,g,h)=(\phi_1(b,c,d,e,g,h),b,c,d,e,g,h).
\]

\paragraph{Step 2. \(v_2=0\) solves \(b=A_{12}\).}
Insert \(E_1(b,c,d,e,g,h)\) into the three products of \(S_2\).
Since \(v_1=0\), the coefficient \(H_{22}\) vanishes, leaving
\[
G_2(b,c,d,e,g,h)=(\Rea\ell S_2)(E_1(b,c,d,e,g,h))=0.
\]
Along \(v_1=0\), implicit differentiation gives
\(\partial_b\phi_1=-v_{1,b}/v_{1,a}\), and hence
\begin{equation}\label{eq:pivot2}
	\partial_bG_2=v_{2,b}-v_{2,a}\frac{v_{1,b}}{v_{1,a}}
	\ \in [-5.0740735,-5.0740734]\times10^{-6}.
\end{equation}
At the certified point, \(v_{2,b}\simeq5.54663\times10^{-4}\), whereas
the derivative after the first substitution is negative. The enclosure in
\eqref{eq:pivot2} gives the local solution \(b=\phi_2(c,d,e,g,h)\) by the
implicit-function theorem. We record this substitution as
\[
E_2(c,d,e,g,h)=E_1(\phi_2(c,d,e,g,h),c,d,e,g,h).
\]

\paragraph{Step 3. \(v_3=0\) solves \(c=A_{13}\).}
Evaluating the six products at \(E_2(c,d,e,g,h)\) gives the next equation
from \(S_3\). This incorporates \(b=\phi_2(c,d,e,g,h)\) together with
\(a=\phi_1(\phi_2(c,d,e,g,h),c,d,e,g,h)\).
The preceding equations \(v_1=v_2=0\) imply \(H_{22}=H_{33}=0\), so we obtain
\[
G_3(c,d,e,g,h)=(\Rea\ell S_3)(E_2(c,d,e,g,h))=0.
\]
Its derivative is
\[
\partial_cG_3
=v_{3,c}-(v_{3,a},v_{3,b})J_2^{-1}(v_{1,c},v_{2,c})^{\mathsf T}
\in[-5.4830652,-5.4830651]\times10^{-7}.
\]
The nonzero derivative yields the local analytic solution
\(c=\phi_3(d,e,g,h)\), with
\(E_3(d,e,g,h)=E_2(\phi_3(d,e,g,h),d,e,g,h)\).

\paragraph{Step 4. \(v_4=0\) solves \(d=A_{21}\).}
With \(E_3(d,e,g,h)\), we substitute \(c\), followed by \(b\) and \(a\).
The first three equations also imply \(H_{44}=0\) in addition to the
earlier diagonal zeros. The ten products in \(S_4\) then give
\(G_4(d,e,g,h)=(\Rea\ell S_4)(E_3(d,e,g,h))=0\).
Differentiating with \(e,g,h\) fixed, we find
\[
\partial_dG_4
=v_{4,d}-(v_{4,a},v_{4,b},v_{4,c})J_3^{-1}
(v_{1,d},v_{2,d},v_{3,d})^{\mathsf T}
\in[5.0676265,5.0676266]\times10^{-10}.
\]
We can therefore solve for \(d=\phi_4(e,g,h)\) and define
\(E_4(e,g,h)=E_3(\phi_4(e,g,h),e,g,h)\).

\paragraph{Step 5. \(v_5=0\) solves \(e=A_{22}\).}
Continuing with \(E_4(e,g,h)\) substitutes \(d,c,b,a\) in that order.
The first four equations imply \(H_{55}=0\). The fifteen products in
\(S_5\) therefore yield \(G_5(e,g,h)=(\Rea\ell S_5)(E_4(e,g,h))=0\).
For fixed \(g,h\), the restricted derivative is
\[
\partial_eG_5
=v_{5,e}-(v_{5,a},v_{5,b},v_{5,c},v_{5,d})J_4^{-1}
(v_{1,e},v_{2,e},v_{3,e},v_{4,e})^{\mathsf T}
\in[-4.1094693,-4.1094692]\times10^{-10}.
\]
It follows that \(e=\phi_5(g,h)\), and the next substitution is
\(E_5(g,h)=E_4(\phi_5(g,h),g,h)\).

\paragraph{Step 6. \(v_6=0\) solves \(g=A_{31}\).}
At this stage, \(E_5(g,h)\) incorporates the solutions for \(e,d,c,b,a\).
The first five equations imply \(H_{66}=0\), leaving the twenty-one
products of \(S_6\). Their sum gives
\(G_6(g,h)=(\Rea\ell S_6)(E_5(g,h))=0\), whose derivative for fixed \(h\) is
\[
\begin{aligned}
	\partial_gG_6
	&=v_{6,g}-(v_{6,a},v_{6,b},v_{6,c},v_{6,d},v_{6,e})J_5^{-1}
	(v_{1,g},v_{2,g},v_{3,g},v_{4,g},v_{5,g})^{\mathsf T}\\
	&\in[-5.3450887,-5.3450886]\times10^{-12}.
\end{aligned}
\]
This determines \(g=\phi_6(h)\) locally and hence
\(E_6(h)=E_5(\phi_6(h),h)\).

\paragraph{Step 7. \(v_7=0\) solves \(h=A_{32}\).}
The final substitution \(E_6(h)\) incorporates \(g,e,d,c,b,a\).
Since the first six equations imply \(H_{77}=0\), the twenty-eight
products in \(S_7\) reduce the last equation to
\(G_7(h)=(\Rea\ell S_7)(E_6(h))=0\). Its derivative satisfies
\[
\begin{aligned}
	G_7'(h)
	&=v_{7,h}-(v_{7,a},v_{7,b},v_{7,c},v_{7,d},v_{7,e},v_{7,g})J_6^{-1}
	(v_{1,h},v_{2,h},v_{3,h},v_{4,h},v_{5,h},v_{6,h})^{\mathsf T}\\
	&\in[5.8000146,5.8000293]\times10^{-13}.
\end{aligned}
\]
Strict positivity makes \(h_*\), the last coordinate of \(p_*\), a simple local root. We have thus
obtained \(p_*=E_6(h_*)\), together with \(f=f_{\rm H}(p_*)\).

\begin{table}[htbp]
	\caption{The seven scalar eliminations and the order-eight sign test. The dyadic interval calculation certifies each displayed derivative along the earlier solution graphs.}
	\label{tab:pivots}
	\centering
	\small
	\renewcommand{\arraystretch}{1.08}
	\setlength{\tabcolsep}{5pt}
	\begin{tabular}{@{}llcl@{}}
		\toprule
		Reduced scalar & Solution of scalar equation & Zero before step & Certified derivative or value\\
		\midrule
		\(G_1(a,b,c,d,e,g,h)\) & \(a=A_{11}=\phi_1(b,c,d,e,g,h)\) & \(H_{11}\) & \([1.0162706,1.0162707]\times10^{-2}\)\\
		\(G_2(b,c,d,e,g,h)\) & \(b=A_{12}=\phi_2(c,d,e,g,h)\) & \(H_{22}\) & \([-5.0740735,-5.0740734]\times10^{-6}\)\\
		\(G_3(c,d,e,g,h)\) & \(c=A_{13}=\phi_3(d,e,g,h)\) & \(H_{33}\) & \([-5.4830652,-5.4830651]\times10^{-7}\)\\
		\(G_4(d,e,g,h)\) & \(d=A_{21}=\phi_4(e,g,h)\) & \(H_{44}\) & \([5.0676265,5.0676266]\times10^{-10}\)\\
		\(G_5(e,g,h)\) & \(e=A_{22}=\phi_5(g,h)\) & \(H_{55}\) & \([-4.1094693,-4.1094692]\times10^{-10}\)\\
		\(G_6(g,h)\) & \(g=A_{31}=\phi_6(h)\) & \(H_{66}\) & \([-5.3450887,-5.3450886]\times10^{-12}\)\\
		\(G_7(h)\) & \(h=A_{32}=h_*\) & \(H_{77}\) & \([5.8000146,5.8000293]\times10^{-13}\)\\
		\(\widehat v_8(p_*)\) & \(\text{Evaluate at }p_*\) & \(H_{88}\) & \(v_8(p_*)>0\)\\
		\bottomrule
	\end{tabular}
\end{table}

\paragraph{Step 8. Evaluate \(v_8\) at the common zero.}
The seven substitutions make \(H_{11},\ldots,H_{88}\) vanish.
We are left with the thirty-six products in \(S_8\).
\[
\begin{aligned}
	v_8(p_*)&=(\Rea\ell S_8)(E_6(h_*))\\
	&\in [1.79627930791107339047187027050,\,
	1.79627930791107339047187027338]\times10^{-16}.
\end{aligned}
\]
At \(p_*\), the first seven focus values vanish and the eighth is
strictly positive. The equilibrium is therefore an eighth-order fine focus.

\subsection{Nonlinearity of the substituted equations}

The Hopf equation \(D=0\) is affine in \(A_{23}\), while the first focus
equation reduces to the degree-seven polynomial in \(A_{11}\) derived
below. For steps 2--7, Table~\ref{tab:curvature} encloses the second
derivative of \(G_j\) with respect to the solved entry \(p_j\), holding
the trailing entries at their critical values. Since every interval
excludes zero at \(p_*\), these \(G_j\) are non-affine in the entries
being solved for. Their local solvability follows from the
first-derivative enclosures in Table~\ref{tab:pivots}.

We compute these derivatives along the solution graphs. Let \(\zeta\)
be the tangent to \(E_{j-1}(p_j,\ldots,p_7)\) when only \(p_j\) varies.
Its first \(j-1\) components are
\(-J_{j-1}^{-1}\partial_{p_j}(v_1,\ldots,v_{j-1})\). The \(j\)-th
component is one, and all remaining components are zero. With
\(Q_i=\nabla^2v_i[\zeta,\zeta]\), a second differentiation of the
preceding zero equations gives the first \(j-1\) acceleration components
as \(-J_{j-1}^{-1}(Q_1,\ldots,Q_{j-1})^{\mathsf T}\). Hence
\[
\partial_{p_j}^2G_j
=Q_j-\frac{\partial v_j}{\partial(p_1,\ldots,p_{j-1})}
J_{j-1}^{-1}(Q_1,\ldots,Q_{j-1})^{\mathsf T}.
\]
For the interval calculation, we differentiate the focus recurrence
automatically on a box centered at \(\mathcal T(p_0)\). Here
\(\mathcal T(p)=p-MF(p)\) is the contraction of Section~\ref{sec:validation},
\(p_0\) is the original box center, \(M\) is its fixed preconditioning
matrix and \(q_B\) bounds its Lipschitz constant in the sup norm.
The estimate
\(\|p_*-\mathcal T(p_0)\|_\infty\le q_B\|MF(p_0)\|_\infty/(1-q_B)\)
follows by applying the contraction bound to \(p_*,p_0\).

\begin{table}[htbp]
	\caption{Second-derivative interval enclosures with the trailing entries fixed at their critical values.}
	\label{tab:curvature}
	\centering
	\small
	\renewcommand{\arraystretch}{1.08}
	\begin{tabular}{@{}ccl@{}}
		\toprule
		Reduced equation & Solved entry & Second derivative at the critical point\\
		\midrule
		\(G_2\) & \(b=A_{12}\) & \([-7.5411741,-7.5411740]\times10^{-5}\)\\
		\(G_3\) & \(c=A_{13}\) & \([3.4246055,3.4246056]\times10^{-5}\)\\
		\(G_4\) & \(d=A_{21}\) & \([-1.5370596,-1.5370595]\times10^{-8}\)\\
		\(G_5\) & \(e=A_{22}\) & \([1.9839111,1.9839112]\times10^{-7}\)\\
		\(G_6\) & \(g=A_{31}\) & \([-3.3764817,-3.3764816]\times10^{-10}\)\\
		\(G_7\) & \(h=A_{32}\) & \([3.07920,3.07924]\times10^{-11}\)\\
		\bottomrule
	\end{tabular}
\end{table}

\subsection{An exact septimic for the first step}\label{sec:septimic}

We now derive the degree-seven equation in \(a=A_{11}\), keeping the
other six entries as parameters. To retain a factored expression for
its coefficients, introduce
\[
N_c=-[(-2\eta-\tau)T_c]_0,\qquad
T_c=\sum_{i,j=1}^3[(A+\eta I)(A-\tau I)]_{ji}(A^2+sI)_{ij}
\mathfrak r_j(\eta)\mathfrak r_i(-\eta),
\]
where \([\cdot]_0\) denotes the constant coefficient after reduction
modulo \(\eta^2+s\). Using \(\mathfrak r(\eta)=Q(\eta)q\) and
\(q\ell=(A+\eta I)(A-\tau I)/[2\eta(\eta-\tau)]\), we obtain
\[
v_1=\frac{N_c}{2(\tau^2+s)(\tau^2+4s)Q(\eta)Q(-\eta)}.
\]
Set
\[
\begin{aligned}
	P_{\rm H}&=\delta f_{\rm H},&
	S_0&=\delta\,s\big|_{f=f_{\rm H}},\\
	L_1&=b(e-1)(a-1)-b^2d-bcg-c(e-1)h,&
	C_3&=\frac{(bP_{\rm H}-ce\delta)^2+c^2\delta S_0}{L_1}.
\end{aligned}
\]
The quantities \(L_1\) and \(C_3\) have degrees one and three in \(a\)
over the rational function field in the remaining six entries.
Polynomial division verifies that \(C_3\) is a polynomial and yields
\[
S_7^{\rm alg}=\frac{\delta^4\,N_c\big|_{f=f_{\rm H}}}{C_3}
\in\mathbb Q[a,b,c,d,e,g,h],\qquad
[a^7]S_7^{\rm alg}=-gh(e-1)^3>0.
\]
Here \([a^7]\) extracts the coefficient of \(a^7\). The degree of \(N_c\) in \(f\)
is four. On the Hopf chart this gives
\begin{equation}\label{eq:septimic}
	v_1=\frac{S_7^{\rm alg}}
	{2L_1(\delta\tau^2+S_0)(\delta\tau^2+4S_0)}.
\end{equation}
Since \(\delta^2Q(\eta)Q(-\eta)=L_1C_3>0\), neither \(L_1\) nor
\(C_3\) vanishes on this chart. The first step is therefore equivalent
to \(S_7^{\rm alg}=0\). At a zero of \(v_1\), the derivative
\(\partial_a S_7^{\rm alg}\) is the denominator in \eqref{eq:septimic}
multiplied by \(\partial_a v_1\). Both factors are nonzero on the
certified branch, so the septimic has a simple root and the
implicit-function theorem selects a unique local branch.

\paragraph{Exact Sturm check on a rational slice.}
For an exact symbolic check, fix \(b,c,d,e,g,h\) at the decimal
rationals in Appendix~\ref{app:parameters} and leave \(a\) free. The polynomial
\(f-f_{\rm H}(a)\) is monic
and linear in \(f\), so its resultant satisfies
\[
\operatorname{Res}_f(f-f_{\rm H}(a),N_1(a,f))=N_1(a,f_{\rm H}(a)).
\]
The two formulas for \(v_1\) give \(N_1=N_cQ(\eta)Q(-\eta)\) on the Hopf
chart, so \(N_1|_{f=f_{\rm H}}=\delta^{-6}L_1C_3^2 S_7^{\rm alg}\).
The specialized numerator and denominator in \eqref{eq:num} have degrees 14 and 12.
Their exact greatest common divisor is, up to a rational constant, a
linear factor times the square of a cubic. After cancellation, the
numerator is irreducible of degree seven over \(\mathbb Q\). The
denominator and \(s\) remain strictly positive on the interval for \(a\)
from the certification box, where Sturm's theorem gives exactly one root of the reduced
numerator. Its monic coefficients, ordered from the leading term to
the constant term and rounded for display, are
\begin{align*}
	(1,&-16.925006124,\ 89.052696790,\ -65.010934159,\ -504.516582752,\\
	&-56.263108911,\ 700.499739173,\ 340.347809648).
\end{align*}
On the same slice, the compact expression has numerator and denominator
degrees 10 and 8, with a cubic common factor and the same reduced
polynomial. This root count concerns the fixed rational slice. The
interval derivatives in Table~\ref{tab:pivots} certify the local
graphs through the simultaneous zero.

\subsection{Center conditions and their scope}\label{sec:center}

On the admissible Hopf surface, a center at \(\mathbf 1\) is
characterized by the simultaneous vanishing of all focus values.
To justify this criterion, choose real coordinates along the center
and stable eigenspaces, and write the center coordinates as
\(r\cos\theta,r\sin\theta\) and the stable coordinate as \(ru\).
The full-turn map \((\mathcal R,\mathcal U)\) in \((r,u)\) is analytic
near \((0,0)\), and \(\partial_u\mathcal U(0,0)=\exp(2\pi\tau/\sqrt{s})\ne1\).
The implicit-function theorem therefore solves \(\mathcal U(r,u)=u\)
analytically for \(u=u(r)\).
A finite-order normal form identifies the leading term of
\(\mathcal R(r,u(r))-r\) at the first nonzero focus value.

For the required order comparison, truncate the center parameterization
at an arbitrary degree \(N\) and combine it with the stable eigenvector
to form a local polynomial coordinate change. The residual of the
stable equation on the center plane is \(O(r^{N+1})\). After blowup,
the stable fixed-point equation gives \(u=O(r^N)\) in these coordinates.
The radial return coefficients through degree \(N\) then agree with those
of the truncated normal form, and changing coordinates preserves the
order of vanishing. If all focus values vanish, the analytic
displacement vanishes to every order and hence identically. This is
the analytic center criterion for the present spectral configuration.
Its relation to local integrability is also treated in \cite{Center3D}.

For the boundary calculations, write \(r_i=-\sum_{j=1}^3 A_{ij}\)
for the intrinsic rates, \(\boldsymbol r=(r_1,r_2,r_3)^{\mathsf T}\), and
\(K_i=-r_i/A_{ii}\) for the density of species \(i\) at its axial
equilibrium. The invasion rate of species \(j\) at \(K_ie_i\) is
\(\lambda_{ij}=r_j+A_{ji}K_i\) for \(i\ne j\). In the class-27 labeling
used here the boundary characteristic number is
\begin{equation}\label{eq:boundaryP}
	P(A)=\lambda_{12}\lambda_{23}\lambda_{31}
	+\lambda_{21}\lambda_{13}\lambda_{32}.
\end{equation}
We write \(P\) when the matrix is fixed.

\begin{proposition}\label{prop:center}
	The certified Hopf box of Section~\ref{sec:validation} contains no center. Each of the following conditions is sufficient for a center at \(E\) in a strictly competitive Hopf system.
	\begin{enumerate}
		\item The intrinsic rates \(r_i\) are equal.
		\item The axial invasion rates satisfy \(\lambda_{ij}+\lambda_{ji}=0\) for every \(i<j\).
	\end{enumerate}
\end{proposition}
\begin{proof}
	The uniform positive enclosure of \(v_8\) excludes centers from the box.
	Under condition 1, set \(u=x_1/x_3,\ v=x_2/x_3\) and use the time
	variable \(\dd\sigma=x_3\dd t\). This gives a planar Lotka--Volterra
	quotient whose positive equilibrium has zero trace and positive
	determinant. The complementary eigenvalue is \(-r_1\), so the planar
	eigenvalues form the Hopf pair. Write the rows of the quotient interaction
	matrix as \((m_{11},m_{12})=(A_{11}-A_{31},A_{12}-A_{32})\) and
	\((m_{21},m_{22})=(A_{21}-A_{31},A_{22}-A_{32})\).
	The multiplier \(u^{\alpha-1}v^{\beta-1}\) makes its divergence vanish
	provided \(m^{\mathsf T}(\alpha,\beta)^{\mathsf T}=-(m_{11},m_{22})^{\mathsf T}\).
	These equations are solvable because \(\det m>0\). The constant
	remaining after the monomial factor is removed is \(\tr m=0\).
	The resulting local Hamiltonian has a nondegenerate extremum.
	Every sufficiently small periodic quotient orbit lifts to a positive
	periodic solution of
	\(\dd x_3/\dd\sigma=r_1+(A_{31}u+A_{32}v+A_{33})x_3\), a scalar linear
	equation with a strictly negative linear coefficient and positive forcing.
	
	For condition 2, put \(y_i=x_i/K_i\) and \(S=\sum_i y_i\), keeping the
	time \(t\). We obtain
	\(\dot y_i=y_i[r_i(1-S)+\sum_j \mathcal M_{ij}y_j]\), with
	\(\mathcal M_{ij}=\lambda_{ji}\) for \(j\ne i\) and \(\mathcal M_{ii}=0\).
	By assumption, \(\mathcal M\) is skew-symmetric, and
	\(\dot S=(1-S)\sum_i r_i y_i\). The positive equilibrium \(y^*\)
	therefore lies on \(S=1\) and satisfies \(\mathcal M y^*=0\).
	On the invariant plane \(S=1\), the first integral
	\(\sum_i y_i^*\log y_i\) has a strict nondegenerate maximum at
	\(y^*\). This gives a center.
\end{proof}

These sufficient families were given by Hofbauer and So \cite{HofbauerSo}.
For example, the circulant matrix with first row
\((-1,-1-\xi,-1+\xi)\), where \(0<\xi<1\), satisfies both conditions.
Its spectrum is \(-3,\ \pm i\sqrt{3}\xi\), and it belongs to class~27.
The full normalized family therefore contains centers.

Each focus value is rational on the admissible Hopf chart. Clear its
nonvanishing denominator and denote the resulting numerator by \(\mathcal N_j\)
in \(\mathbb Q[a,b,c,d,e,g,h]\). The Hilbert basis theorem implies
that the increasing sequence \((\mathcal N_1,\ldots,\mathcal N_m)\) eventually equals
the ideal generated by all focus-value numerators. A finite initial
segment of the focus values therefore characterizes the center set
on this chart. At \(p_*\), we have \(\mathcal N_1(p_*)=\cdots=\mathcal N_7(p_*)=0\)
and \(\mathcal N_8(p_*)\ne0\). A test based only on an initial segment of focus
values therefore requires at least eight. This argument leaves open the
stabilization index and the full real class-27 center locus. For center varieties of generalized three-dimensional
Lotka--Volterra systems, see \cite{BobZol}.

A class-27 center must also satisfy \(P=0\). To see this, suppose
\(P\ne0\). The weighted-product construction of Section~\ref{sec:exterior} gives
a collar containing no positive periodic orbit. When \(P<0\), we
reverse all three weight inequalities, so the averaged logarithmic
derivative has the opposite strict sign. A periodic orbit entering a
sufficiently small sublevel in the collar would attain a positive
minimum there, contradicting that derivative.

If \(E\) were a center, extend its neighborhood through the neighboring
open periodic annuli to obtain the connected open region containing
\(E\) in the carrying simplex. Its closure remains in the interior
of the carrying simplex. The boundary is nonempty and invariant,
and contains no equilibrium, since \(E\) is unique and lies in the
interior of the region. Poincar\'e--Bendixson therefore gives a periodic
orbit on this boundary. Periodic trajectories from the region
accumulate on that orbit, giving accumulating one-turn fixed points
on a short transversal as in Section~\ref{sec:exterior}. The analytic return-map
alternative there gives an open periodic neighborhood of the boundary orbit. This
neighborhood meets and extends the region, contradicting maximality.

The diagonal of \(A\operatorname{diag}(K_1,K_2,K_3)+\boldsymbol r\mathbf1^{\mathsf T}\)
is zero, and its off-diagonal entries are \(\lambda_{ji}\). Its determinant
is therefore \(P\). Since \(A^{-1}\boldsymbol r=-\mathbf1\), the matrix
determinant lemma yields
\[
P=\det A\,K_1K_2K_3\left(1-\sum_{i=1}^3 K_i^{-1}\right).
\]
Consequently, a class-27 center requires
\(\sum_i(-A_{ii})/r_i=1\). Even with \(P=0\), vanishing of the first
two focus values on the Hopf surface does not ensure a center
\cite[Theorem~2.1]{GYDEA}.

\section{Verification of the critical point and the unfolding condition}\label{sec:validation}

The critical matrix is \(A_*=A(p_*)\), where \(p_*\) is the zero of \(F\)
certified below. Its entries rounded to 12 decimal places are
\[
A_*=\begin{pmatrix}
	-0.579177964532 & -2.934485224627 & -1.321825219182\\
	-0.725374736374 & -3.820698848362 & -1.025598329826\\
	-0.395539625144 & -3.746192929039 & -1
\end{pmatrix}.
\]
These digits locate the example. With the seven free entries rounded as shown
and \(A_{23}\) recomputed from \eqref{eq:hopf}, the first focus value is already
about \(8.11\times10^{-16}\). Let \(p_0\) be the exact rational center
given in Appendix~\ref{app:parameters}, and set
\(B=\{p:\|p-p_0\|_\infty\le r_B\}\), where \(r_B=2^{-132}\) and
\(\|\cdot\|_\infty\) is the sup norm, with its induced matrix norm.

\begin{lemma}\label{lem:cert}
	There is a unique \(p_*\) in \(B\) with \(F(p_*)=0\). The Jacobian \(DF\)
	on the Hopf chart is invertible throughout \(B\), and \(v_8(p_*)>0\).
\end{lemma}
\begin{proof}
	All interval endpoints are rounded outward to multiples of \(2^{-384}\).
	Evaluation of \eqref{eq:R}--\eqref{eq:H} and differentiation with
	respect to \(p\) enclose \(DF(p)\) for \(p\in B\). We obtain the fixed
	dyadic matrix \(M\) by rounding the inverse of the midpoint matrix of
	this enclosure. The bounds are
	\[
	\sup_{p\in B}\|I-MDF(p)\|_\infty<q_B=4.644\times10^{-17},
	\qquad \|MF(p_0)\|_\infty<e_B=1.020926682\times10^{-63}.
	\]
	The map \(\mathcal T(p)=p-MF(p)\) moves the box center by at most \(e_B\)
	and has Lipschitz constant at most \(q_B\) on the box. Since
	\(q_Br_B+e_B<r_B\), it sends the box strictly into itself. The Neumann
	series makes \(MDF(p)\) invertible throughout \(B\), so both \(M\) and
	\(DF(p)\) are invertible. The unique fixed point is therefore a zero of
	\(F\), and satisfies
	\[
	\|p_*-p_0\|_\infty\le \frac{e_B}{1-q_B}<1.020926683\times10^{-63}.
	\]
	On the original box the interval for the eighth focus value is
	\(v_8(B)\subset[1.79602,1.79654]\times10^{-16}\). Re-evaluating on the
	tighter box supplied by the fixed-point bound gives
	\[
	v_8(p_*)\in
	[1.79627930791107339047187027050,\,
	1.79627930791107339047187027338]\times10^{-16}.
	\]
\end{proof}

Multiplying the dyadic pivot intervals gives
\[
\det DF(p)\in[1.8254244,1.8254291]\times10^{-56}\qquad(p\in B).
\]
All nine matrix entries are strictly negative. At the critical point
\(\tau\approx-5.39988\) and \(s\approx0.119216\).

To vary the real part of the conjugate eigenvalues, use
\(\xi=f-f_{\rm H}(p)\) as the parameter transverse to the Hopf chart.
Near \(A_*\), write the eigenvalues as
\(\gamma,\alpha\pm i\omega\), with \(\gamma<0,\ \omega>0\), and put
\(v_0(\xi,p)=\alpha\).

Adjoin \((\xi,p)\) as variables with zero time derivatives and take a
local \(C^{20}\) center manifold for this extended system. Each slice
at fixed \((\xi,p)\) is two-dimensional. On it we use the finite normal
form
\begin{equation}\label{eq:finite-normal-form}
	\dot z=(\alpha+i\omega)z+\sum_{j=1}^8c_j(\xi,p)z|z|^{2j}
	+\mathcal R_{\rm nf}(z,\bar z,\xi,p),
	\qquad \mathcal R_{\rm nf}=O(|z|^{18}),\quad
	D_z\mathcal R_{\rm nf}=O(|z|^{17}).
\end{equation}
Here \(D_z\) denotes the derivative in the two real coordinates of \(z\),
and the estimates are uniform for nearby parameters. We retain
\(z^{j+1}\bar z^j\) as in \eqref{eq:H}, even off the Hopf chart. The
denominators for the other monomials remain nonzero near \(A_*\).
Thus \(v_j(\xi,p)=\Rea c_j(\xi,p)\) for \(1\le j\le8\) are smooth
extensions of the focus values on \(\xi=0\).

The identity \(D=-2\alpha[(\gamma+\alpha)^2+\omega^2]\) gives
\[
\left.\frac{\partial v_0}{\partial A_{23}}\right|_{\mathrm{Hopf}}
=-\frac{\delta}{2(\tau^2+s)}\in[-0.32823363,-0.32823362]
\qquad(p\in B).
\]
On \(\xi=0\), the identity \(v_0(0,p)=0\) makes every \(p\)-entry of
the first Jacobian row vanish, and the lower right block is \(DF\).
Since the change from \((f,p)\) to \((\xi,p)\) has determinant one,
the original column order \((f,a,b,c,d,e,g,h)\) gives
\[
\det\frac{\partial(v_0,v_1,\ldots,v_7)}{\partial(f,a,b,c,d,e,g,h)}
=-\frac{\delta}{2(\tau^2+s)}\det DF
\in[-5.9916720,-5.9916567]\times10^{-57}.
\]
The enclosure uses the unrounded interval factors. In particular, the
inverse function theorem gives a smooth local parameterization by
\((v_0,\ldots,v_7)\) near the critical values.

\section{Hopf cyclicity and the exterior cycle}\label{sec:cycles}

\subsection{Eight small-amplitude limit cycles}

A sufficiently small periodic orbit stays in the center-manifold
neighborhood for all positive and negative times and hence lies in
the manifold. We use \eqref{eq:finite-normal-form} on the slice for
each fixed parameter and write \(z=re^{i\theta}\).

Let \(\omega_*=\sqrt{s(p_*)}\) and \(\nu=v_8(p_*)>0\).
Write \(\mathsf e_m\) for the \(m\)th elementary symmetric polynomial.
For small \(\varepsilon>0\), use the local inverse from
Section~\ref{sec:validation} to choose
\[
v_j=\nu(-1)^{8-j}\mathsf e_{8-j}(1,2,\ldots,8)\,\varepsilon^{2(8-j)},
\qquad 0\le j\le7,
\]
so the parameters change by \(O(\varepsilon^2)\). These coefficients give the
eight prescribed roots in the leading radial polynomial. Indeed,
substituting \(r=\varepsilon\rho\) and dividing the radial equation by
the angular speed gives
\[
\frac{\dd\rho}{\dd\theta}
=\varepsilon^{16}\frac{\nu}{\omega_*}
\rho\prod_{m=1}^8(\rho^2-m)+o_{C^1}(\varepsilon^{16})
\]
where the \(C^1\) error is in \(\rho\), uniformly in \(\theta\) on compact
subintervals of \(\rho>0\). Here
the radial terms through degree 17 are \(v_jr^{2j+1}\),
\(v_8=\nu+O(\varepsilon^2)\), and the angular speed at this scale is
\(\omega_*+O(\varepsilon^2)\).
Let \(\Pi_\varepsilon(r)\) be the radial return after one turn for
these parameters. Integrating the scaled equation and converting
back to \(r=\varepsilon\rho\) gives, in \(C^1\) with respect to \(\rho\),
\[
\varepsilon^{-17}\big(\Pi_\varepsilon(\varepsilon\rho)-\varepsilon\rho\big)
\longrightarrow \frac{2\pi\nu}{\omega_*}\rho\prod_{m=1}^8(\rho^2-m).
\]
The eight positive simple zeros give eight hyperbolic small-amplitude
limit cycles. At the cycle indexed by \(m\), the radial multiplier
minus one has the sign of \(2m\prod_{\ell\ne m}(m-\ell)\). In
particular the outermost small cycle repels radially. Its other
nontrivial multiplier tends to \(\exp(2\pi\tau(p_*)/\omega_*)<1\).

The full-turn map also gives a local upper bound. Use a signed radius
on a line through the equilibrium in the slice, with angular speed
bounded away from zero on a uniform neighborhood. The scalar map
preserves the two half-lines and is increasing on each one, by
uniqueness for the differential equation with angle as independent
variable. Consequently a periodic point of the map is fixed, and each
small periodic orbit contributes one fixed point on each half-line.
At the critical parameters the signed displacement is
\[
\frac{2\pi\nu}{\omega_*}r^{17}+o(r^{17}).
\]
Its continuous seventeenth derivative in \(r\) takes the value
\(17!\,2\pi\nu/\omega_*\ne0\) at the critical parameters and \(r=0\).
It remains nonzero on a smaller radius and parameter neighborhood.
Rolle's theorem allows at most 17 distinct signed zeros, one of which is
the equilibrium. The eight cycles just constructed attain the resulting
bound, giving local cyclicity exactly eight at the critical equilibrium.

\subsection{A ninth cycle in the same perturbed system}\label{sec:exterior}

For the invasion rates \(\lambda_{ij}\) defined in
Section~\ref{sec:center}, the certified signs in the order
\(12,13,21,23,31,32\) are \((-,+,+,-,-,+)\). They identify the class-27
boundary cycle in our labeling. The interval certificate at the
critical parameters gives \(P(A_*)\in[5.51998,5.52000]\times10^{-10}\)
for the characteristic number in \eqref{eq:boundaryP}.
The invasion signs and \(P>0\) persist for the perturbation used to
produce the eight small cycles. By the boundary repulsion criterion,
the heteroclinic cycle repels trajectories in the carrying simplex
\cite{GYDEA,JiangNiu}.

\paragraph{A compact region away from the boundary.}
Choose weights \(w_i>0\) such that
\[
\frac{w_3}{w_2}>\frac{-\lambda_{12}}{\lambda_{13}},\qquad
\frac{w_1}{w_3}>\frac{-\lambda_{23}}{\lambda_{21}},\qquad
\frac{w_2}{w_1}>\frac{-\lambda_{31}}{\lambda_{32}} .
\]
Since \(P>0\), the product of the right sides is less than one,
which permits the three strict inequalities. Put
\(\mathcal V(x)=\prod_i x_i^{w_i}\) and
\(\psi(x)=\sum_i w_i[A(x-\mathbf1)]_i\). Then
\(\dot{\mathcal V}/\mathcal V=\psi\) in the positive octant, and the weight inequalities
give \(\psi>0\) at the axial equilibria.
Let \(\varphi_u\) denote the flow for the chosen parameters. Every
invariant probability measure on the boundary of the carrying
simplex is supported on the three saddles, since the connecting
orbits are nonrecurrent. Thus for some common time \(T_0>0\),
\(T_0^{-1}\int_0^{T_0}\psi(\varphi_u x)\,\dd u>0\) at every boundary point.
If such a time did not exist, time averages based at boundary points
over intervals tending to infinity would have an invariant weak
limit with nonpositive integral of \(\psi\). The positive average also
holds in a collar by continuity.

Put
\[
\widetilde{\mathcal V}(x)=\mathcal V(x)\exp\left(
\frac1{T_0}\int_0^{T_0}(T_0-u)\psi(\varphi_u x)\,\dd u\right).
\]
The logarithmic derivative of \(\widetilde{\mathcal V}\) along
interior trajectories is the time average just obtained.
The function \(\widetilde{\mathcal V}\) extends continuously to zero on the
boundary. Each sufficiently small positive level lies in the collar,
where trajectories cross it only in the increasing direction.
An interior trajectory starting above this level therefore stays in
a compact subset of the interior.

\paragraph{The exterior limit set.}
For a positive periodic orbit \(x(t)\), take \(y_i=\log x_i\).
The forward variational matrix in these coordinates is
\(A\,\operatorname{diag}(x(t))\). Every entry changes to a positive
one in reverse time. The reverse monodromy thus has positive entries
and diagonal entries greater than one. Its simple Perron multiplier
exceeds one, so the forward monodromy has a simple multiplier in
\((0,1)\) with a positive eigenvector at each phase. This positivity
persists in the original coordinates because the derivative of the
inverse logarithm is a positive diagonal matrix.

For the outermost small cycle, the other transverse multiplier exceeds
one, so the positive line just found is its stable direction. The local
stable manifold is foliated by one-dimensional stable fibers.
Sufficiently near the cycle, an off-cycle point of a fiber is strictly
coordinatewise comparable to its phase point on the cycle.
Unorderedness of the carrying simplex excludes all such points.
In particular, a trajectory in the carrying simplex and outside the
outermost small cycle cannot approach that cycle.

Radial projection identifies the carrying simplex with a planar
triangle. In the interior the projected vector field is locally
Lipschitz, since the inverse projection is locally Lipschitz and the
ambient field and the projection are smooth. Take a short transverse
arc with a common crossing orientation, and choose a point just
outside the outermost cycle whose one-turn return is on the arc and
differs from the starting point. Hyperbolicity excludes other fixed
points close to the cycle. Successive crossings of a periodic Jordan
curve along the arc would have to alternate directions. The flow crosses
in the same direction throughout the arc, so a periodic orbit crosses
it at most once. The trajectory we have chosen is therefore nonperiodic.

The trajectory cannot cross the outermost small cycle and stays away
from the boundary by the preceding level estimate. The only positive
equilibrium is the one determined by the nonsingular \(A\). Each
small cycle encloses it on the carrying simplex by the planar index
theorem. The exterior limit set thus contains no equilibrium.
Poincar\'e--Bendixson makes it a periodic orbit, and the stable-fiber
argument places it strictly outside all eight small cycles.

\paragraph{Isolation of the exterior orbit.}
Take an analytic transverse section in the full flow. The multiplier
in \((0,1)\) lets us solve the stable component of the fixed point
equation by the implicit-function theorem. The remaining scalar
displacement is analytic. If it vanished identically, the local
fixed curve would sweep out an open annulus of periodic trajectories
in the carrying simplex, by invariance of domain. The approaching
trajectory is nonperiodic and could not enter this annulus.
Consequently the displacement has an isolated zero. This is the
local alternative used by Xiao and Li \cite[proof of Theorem~2.6]{XiaoLi}.
On a sufficiently small relative transverse arc of the carrying
simplex, a periodic orbit starting close enough returns first to that
arc. As above, the common crossing orientation and the Jordan curve
theorem require its return point to be its starting point. The
exterior periodic orbit is therefore a limit cycle, proving
Theorem~\ref{thm:main}.

Fix a sufficiently small perturbation from the construction.
The eight hyperbolic cycles persist under further small changes of
the matrix, with their nesting and multiplier inequalities intact.
The strict invasion signs and \(P>0\) are also open conditions, so
the exterior argument gives the ninth cycle throughout a neighborhood.
The lower bound of nine thus holds on a nonempty open set of the
normalized family. Density then gives rational matrices with at least
nine cycles.

\subsection{Numerical illustration}

Figures~\ref{fig:orbits} and~\ref{fig:orbits3d} show numerical trajectories
for the rational matrix \(A_{\rm num}\) in
Appendix~\ref{app:parameters}. Exact rational arithmetic verifies
strict competition, the class-27 invasion signs and
\(P(A_{\rm num})\in[5.52023,5.52024]\times10^{-10}\).
We selected the matrix using a truncated radial normal form, then
used Poincar\'e shooting in the original three-dimensional system
to compute an initial condition and period for each plotted orbit.

Near \(E\), we solved the return equation for the transverse section
coordinate first. The resulting scalar displacement alternated in
sign at nine ordered section points, giving eight numerical root
brackets. The signs agreed in 340-bit Taylor integrations of
degrees 38 and 48 with respective step sizes \(0.19\) and \(0.11\).
Over the nine points, the smallest absolute displacement was about
\(1.28\times10^{-52}\) and the largest observed difference between
runs was about \(2.11\times10^{-74}\).

For the exterior orbit, shooting on \(x_1=x_3\) gave a period of
approximately \(22.51670160585\). We checked it by halving the step
in logarithmic coordinates and repeating the computation in polynomial
coordinates. We integrated the plotted trajectories in double precision
over the resulting periods. The existence assertion of
Theorem~\ref{thm:main} follows from the interval verification of
the fine focus and the perturbation proof in
Sections~\ref{sec:validation} and~\ref{sec:cycles}.

\begin{figure}[!htbp]
	\centering
	\includegraphics[width=\linewidth]{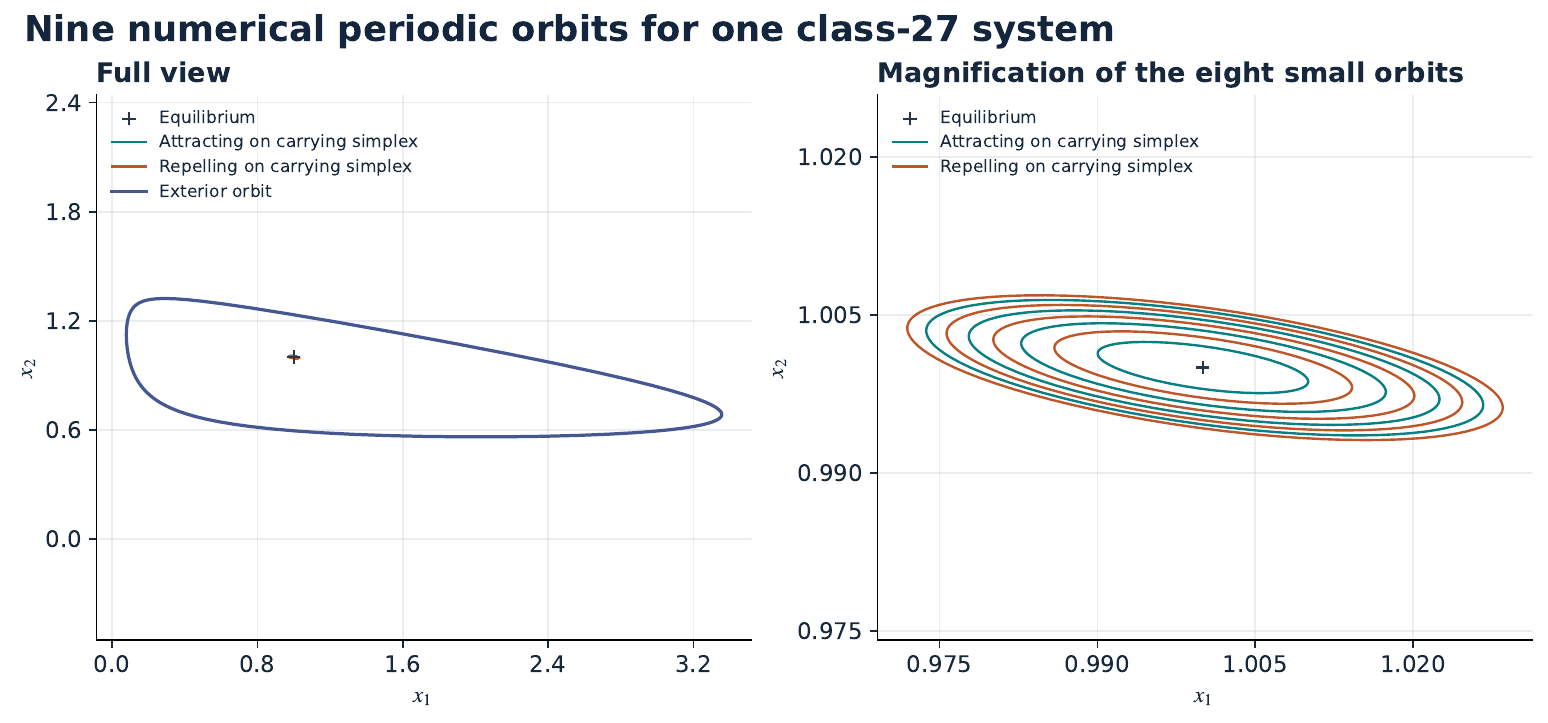}
	\caption{Projection of the numerical trajectories onto the \((x_1,x_2)\)-plane. The left panel shows the exterior orbit (blue) and the eight small-amplitude orbits near \(E=(1,1,1)\), marked by a cross. The right panel magnifies the small-amplitude orbits. Teal and orange indicate numerically attracting and repelling small-amplitude orbits on the carrying simplex, respectively.}
	\label{fig:orbits}
\end{figure}

Figure~\ref{fig:orbits3d} shows the same trajectories in all three
species coordinates, with equal units on the three axes of each panel.
A second integration using a smaller maximum step and tighter
tolerances gave a maximum sampled coordinate difference of
\(2.18\times10^{-10}\). The two integrations agree at the plotting scale.

\begin{figure}[!htbp]
	\centering
	\includegraphics[width=\linewidth]{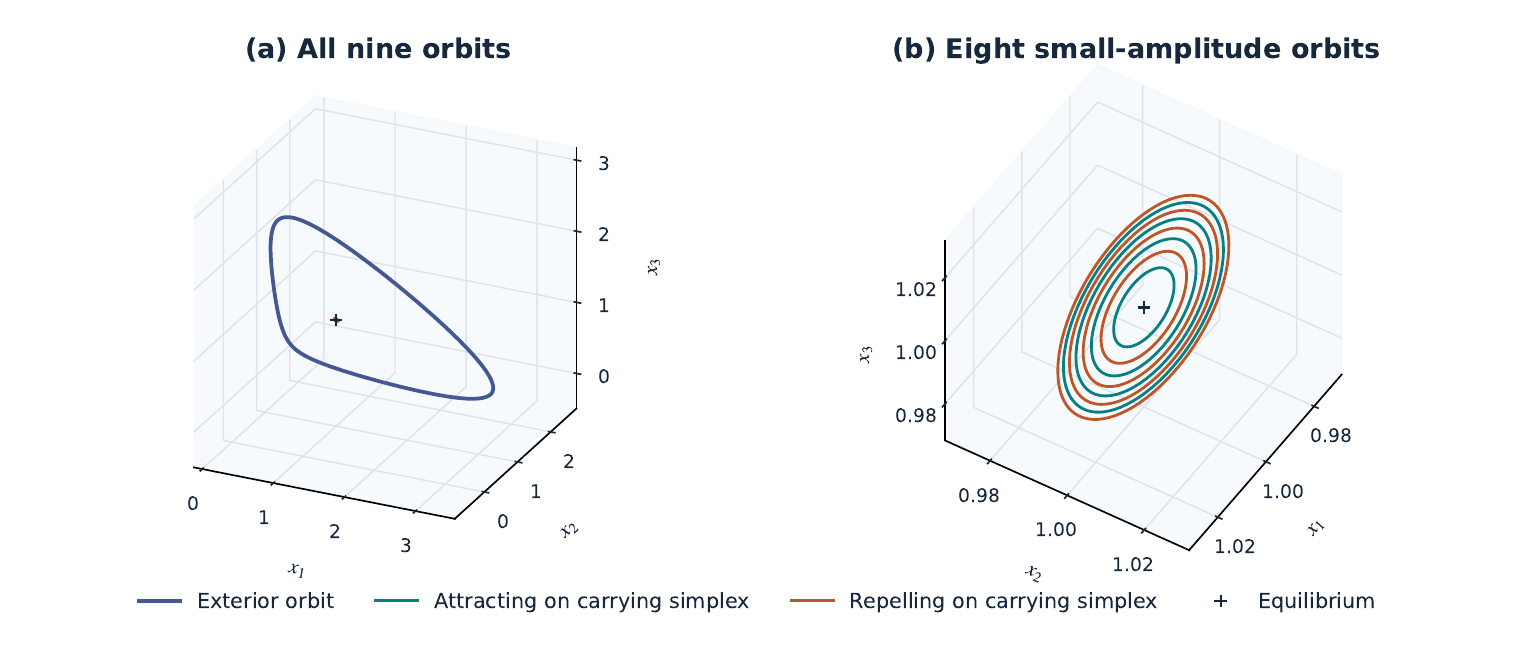}
	\caption{Three-dimensional phase portrait of the same numerical orbits as in Fig.~\ref{fig:orbits}. Each curve shows \((x_1(t),x_2(t),x_3(t))\) over one numerically computed period. Panel (a) shows the exterior orbit and the eight small-amplitude orbits, which are tightly clustered near \(E\) on this scale. Panel (b) magnifies the eight small-amplitude orbits from a different viewing angle. The cross and the colors have the same meaning as in Fig.~\ref{fig:orbits}.}
	\label{fig:orbits3d}
\end{figure}

\section{Conclusion}

The critical matrix gives an eighth-order fine focus, and the unfolding realizes its full local cyclicity of eight. Under the same perturbation, the outermost of the eight small cycles repels on the carrying simplex. The boundary heteroclinic cycle repels as well. These two signs give an enclosing cycle, so the resulting class~27 system has at least nine limit cycles.

A center in class~27 requires \(P=0\) and the vanishing of all focus values. The example also shows that a center test based on the focus values alone must go at least as far as the eighth value. The stabilization index of the ideal generated by the focus-value numerators remains open.

\appendix
\section{Terms retained in the eight focus values}\label{app:terms}
We write \(ab\mid cd\) for \(\BB(H_{ab},H_{cd})\), with each digit denoting one index. The first column gives the lower degree \(m\), and the products on that row make up \(T_{jm}\). Adding the rows at order \(j\) gives \(S_j\) in Table~\ref{tab:focusforms}. For the full value, include the additional terms in \eqref{eq:reduced}. The entries below number \(1,3,6,10,15,21,28,36\) at the successive orders.
\par\smallskip\noindent\textbf{Order 1.}\par\nopagebreak
\begin{tabular}{@{}r l@{}}
	1&$01\mid20$ \\
\end{tabular}
\par\smallskip\noindent\textbf{Order 2.}\par\nopagebreak
\begin{tabular}{@{}r l@{}}
	1&$01\mid31$ \\
	2&$02\mid30,\ 20\mid12$ \\
\end{tabular}
\par\smallskip\noindent\textbf{Order 3.}\par\nopagebreak
\begin{tabular}{@{}r l@{}}
	1&$01\mid42$ \\
	2&$02\mid41,\ 20\mid23$ \\
	3&$03\mid40,\ 12\mid31,\ 30\mid13$ \\
\end{tabular}
\par\smallskip\noindent\textbf{Order 4.}\par\nopagebreak
\begin{tabular}{@{}r l@{}}
	1&$01\mid53$ \\
	2&$02\mid52,\ 20\mid34$ \\
	3&$03\mid51,\ 12\mid42,\ 30\mid24$ \\
	4&$04\mid50,\ 13\mid41,\ 31\mid23,\ 40\mid14$ \\
\end{tabular}
\par\smallskip\noindent\textbf{Order 5.}\par\nopagebreak
\begin{tabular}{@{}r l@{}}
	1&$01\mid64$ \\
	2&$02\mid63,\ 20\mid45$ \\
	3&$03\mid62,\ 12\mid53,\ 30\mid35$ \\
	4&$04\mid61,\ 13\mid52,\ 31\mid34,\ 40\mid25$ \\
	5&$05\mid60,\ 14\mid51,\ 23\mid42,\ 41\mid24,\ 50\mid15$ \\
\end{tabular}
\par\smallskip\noindent\textbf{Order 6.}\par\nopagebreak
\begin{tabular}{@{}r l@{}}
	1&$01\mid75$ \\
	2&$02\mid74,\ 20\mid56$ \\
	3&$03\mid73,\ 12\mid64,\ 30\mid46$ \\
	4&$04\mid72,\ 13\mid63,\ 31\mid45,\ 40\mid36$ \\
	5&$05\mid71,\ 14\mid62,\ 23\mid53,\ 41\mid35,\ 50\mid26$ \\
	6&$06\mid70,\ 15\mid61,\ 24\mid52,\ 42\mid34,\ 51\mid25,\ 60\mid16$ \\
\end{tabular}
\par\smallskip\noindent\textbf{Order 7.}\par\nopagebreak
\begin{tabular}{@{}r l@{}}
	1&$01\mid86$ \\
	2&$02\mid85,\ 20\mid67$ \\
	3&$03\mid84,\ 12\mid75,\ 30\mid57$ \\
	4&$04\mid83,\ 13\mid74,\ 31\mid56,\ 40\mid47$ \\
	5&$05\mid82,\ 14\mid73,\ 23\mid64,\ 41\mid46,\ 50\mid37$ \\
	6&$06\mid81,\ 15\mid72,\ 24\mid63,\ 42\mid45,\ 51\mid36,\ 60\mid27$ \\
	7&$07\mid80,\ 16\mid71,\ 25\mid62,\ 34\mid53,\ 52\mid35,\ 61\mid26,\ 70\mid17$ \\
\end{tabular}
\par\smallskip\noindent\textbf{Order 8.}\par\nopagebreak
\begin{tabular}{@{}r l@{}}
	1&$01\mid97$ \\
	2&$02\mid96,\ 20\mid78$ \\
	3&$03\mid95,\ 12\mid86,\ 30\mid68$ \\
	4&$04\mid94,\ 13\mid85,\ 31\mid67,\ 40\mid58$ \\
	5&$05\mid93,\ 14\mid84,\ 23\mid75,\ 41\mid57,\ 50\mid48$ \\
	6&$06\mid92,\ 15\mid83,\ 24\mid74,\ 42\mid56,\ 51\mid47,\ 60\mid38$ \\
	7&$07\mid91,\ 16\mid82,\ 25\mid73,\ 34\mid64,\ 52\mid46,\ 61\mid37,\ 70\mid28$ \\
	8&$08\mid90,\ 17\mid81,\ 26\mid72,\ 35\mid63,\ 53\mid45,\ 62\mid36,\ 71\mid27,\ 80\mid18$ \\
\end{tabular}

\section{Exact parameter values}\label{app:parameters}
The entries below are the exact coordinates of the box center \(p_0\).
In each coordinate the radius is \(r_B=2^{-132}\). The full Hopf matrix
at each point of the box is obtained with \(k=-1\) and \(f=f_{\rm H}(p)\).
\begin{center}
	\small
	\renewcommand{\arraystretch}{1.08}
	\begin{tabular}{@{}ll@{}}
		\toprule
		Entry & Exact coordinate of the box center \\
		\midrule
		\(A_{11}\) & \texttt{-0.5791779645319467551172476131355811359222106253213525266771896548} \\
		\(A_{12}\) & \texttt{-2.934485224626673940708888971277265018886343541896677486144401152} \\
		\(A_{13}\) & \texttt{-1.321825219182311398189873585735460420576694184584086937916485516} \\
		\(A_{21}\) & \texttt{-0.7253747363735585867032134309525727120024623667097688021323810841} \\
		\(A_{22}\) & \texttt{-3.820698848362104024804343978614026348544247656786447940223496871} \\
		\(A_{31}\) & \texttt{-0.3955396251443007826543881818804962477423972983010672235309116071} \\
		\(A_{32}\) & \texttt{-3.746192929038915742803317816415078247032099360498182764013265698} \\
		\bottomrule
	\end{tabular}
\end{center}

\par\medskip\noindent
\begin{minipage}{\linewidth}
	\textbf{Matrix used for the figures.}
	We set the \((3,3)\) entry of \(A_{\rm num}\) to \(-1\) and give the other entries as exact decimals.
	\begin{center}
		\small
		\renewcommand{\arraystretch}{1.08}
		\begin{tabular}{@{}ll@{}}
			\toprule
			Entry & Exact value \\
			\midrule
			\(A_{11}\) & \texttt{-0.57917810395343374357314547283441548815614384268441597627492438341889} \\
			\(A_{12}\) & \texttt{-2.9344862001292182581951345824963109906094686127607099718229032286987} \\
			\(A_{13}\) & \texttt{-1.3218249495537193478694907393589751268236620003396844681477586394354} \\
			\(A_{21}\) & \texttt{-0.72537473546865998768978328085201908506848331337059392753178433625916} \\
			\(A_{22}\) & \texttt{-3.820698882805152456422613798118902445840929390314805729726085822463} \\
			\(A_{23}\) & \texttt{-1.0255982493303249385344307794975782494243699978863796467054641568241} \\
			\(A_{31}\) & \texttt{-0.39553987156112489629130184816311704231538128567436660164686538380392} \\
			\(A_{32}\) & \texttt{-3.7461932075610644153752496182987881896348018217897762309374215175055} \\
			\bottomrule
		\end{tabular}
	\end{center}
\end{minipage}
\par

\end{document}